\documentclass[12pt,a4paper]{amsart}
\usepackage{a4wide}
\usepackage[english]{babel}
\usepackage[T2A]{fontenc}
\usepackage[utf8]{inputenc} 
\usepackage{amsfonts}
\usepackage{amssymb, amsthm, amscd}
\usepackage{amsmath}
\usepackage{mathtools}
\usepackage{needspace}
\usepackage{etoolbox}
\usepackage{lipsum}
\usepackage{comment}
\usepackage{cmap}
\usepackage[pdftex]{graphicx}
\usepackage[unicode]{hyperref}
\usepackage[matrix,arrow,curve]{xy}
\usepackage[usenames,dvipsnames]{xcolor}
\usepackage{colortbl}
\usepackage{textcomp}
\usepackage{cite}
\usepackage{euscript}
\usepackage{epigraph}
\usepackage[numbers]{natbib}

\newcommand{\Z}{\mathbb{Z}}

\newcommand{\N}{\mathbb{N}}

\newcommand{\Q}{\mathbb{Q}}

\newcommand{\F}{\mathbb{F}}

\newcommand{\ph}{\varphi}

\newcommand{\sub}{\subseteq}
\newcommand{\Hom}{\mathop{\mathrm{Hom}}\nolimits}
\newcommand{\Ker}{\mathop{\mathrm{Ker}}\nolimits}

\newcommand{\rk}{\mathop{\mathrm{rk}}\nolimits}

\renewcommand{\ge}{\geqslant}
\renewcommand{\le}{\leqslant}
\newcommand{\sm}{\setminus}
\newcommand{\map}[3]{#1\colon #2\to #3}

\renewcommand{\>}{\rangle}

\renewcommand{\P}{\mathfrak{p}}

\newcommand{\SL}{\mathop{\mathrm{SL}}\nolimits}
\newcommand{\GL}{\mathop{\mathrm{GL}}\nolimits}

\newcommand{\tc}{\text{,}}
\newcommand{\tp}{\text{.}}

\theoremstyle{plain}
\newtheorem{thm}{Theorem}[section]
\newtheorem{lem}[thm]{Lemma}

\theoremstyle{definition}

\theoremstyle{remark}

\AtBeginEnvironment{thm}{\begin{samepage}}
	\AtEndEnvironment{thm}{\end{samepage}}
\AtBeginEnvironment{lem}{\begin{samepage}}
	\AtEndEnvironment{lem}{\end{samepage}}
\AtBeginEnvironment{st}{\begin{samepage}}
	\AtEndEnvironment{st}{\end{samepage}}
\AtBeginEnvironment{crit}{\begin{samepage}}
	\AtEndEnvironment{crit}{\end{samepage}}
\AtBeginEnvironment{ax}{\begin{samepage}}
	\AtEndEnvironment{ax}{\end{samepage}}
\AtBeginEnvironment{defn}{\begin{samepage}}
	\AtEndEnvironment{defn}{\end{samepage}}
\AtBeginEnvironment{cor}{\begin{samepage}}
	\AtEndEnvironment{cor}{\end{samepage}}
\AtBeginEnvironment{note}{\begin{samepage}}
	\AtEndEnvironment{note}{\end{samepage}}
\AtBeginEnvironment{prop}{\begin{samepage}}
	\AtEndEnvironment{prop}{\end{samepage}}

\makeatletter
\def\@settitle{\begin{center}%
		\baselineskip14\p@\relax
		\bfseries
		\@title
	\end{center}%
}

\def\@evenhead{\hfil\sc pavel gvozdevsky\hfil}
\def\@oddhead{\hfil\sc verbal width\hfil}
\makeatother

\title{Verbal width in arithmetic Chevalley groups}


\author{Pavel Gvozdevsky*}
\date{}
\address{Department of Mathematics, Bar-Ilan University, 5290002 Ramat Gan, ISRAEL}
\thanks{The paper is written as part of the author's post-doctoral fellowship at Bar-Ilan University, Department of Mathematics; and is supported by ISF grant 1994/20.}
\email{gvozdevskiy96@gmail.com}

\begin{document}

\maketitle

\epigraph{Dedicated to the memory of Nikolai Vavilov}

\begin{abstract}
We prove that the width of any word in a simply connected Chevalley group of rank at least 2 over the ring that is a localisation of the ring of integers in a number field is bounded by a constant that depends only on the root system and on the degree of the number field. 
\end{abstract}

\section{Introduction}

In the remarkable paper~\cite{AvniMeiri}, Nir Avni and Chen Meiri studied the width of words in $\SL(n,\Z)$. In particular, they proved \cite[Theorem~1.2]{AvniMeiri} that the width of any word in $\SL(n,\Z)$, $n\ge 3$ is bounded by some constant that depend only on $n$. In the present paper, being inspired by this result, we extend it in the following ways. Firstly, we generalise it to all Chevalley groups. Secondly, we generalize it to all arithmetic rings. Thirdly, we simplify some parts of the proof. Now we explain the setting of the problem in details.

\smallskip   

For a group $G$ {\it a word map} is a map
$$
G^n\to G\tc\qquad
(g_1,\ldots,g_m)\mapsto w(g_1,\ldots,g_m)\tc
$$     
where $w$ is a word of $m$ variables (i.e an element of the free group on $m$ generators), and where $w(g_1,\ldots,g_n)$ is what we get if we substitute the elements $g_1$,$\ldots$,$g_n$ into this word.

In recent studies of word maps, the case of $G=\mathcal{G}(K)$, where $\mathcal{G}$ is a simple or semisimple algebraic group defined over a field $K$, is of particular interest, see references in \cite{GordKunPlotJIJAC}, \cite{GordKunPlotJAlg},\cite{GordKunPlotRAN}. Almost all the results in this direction have been proved for the cases where $K$ is an algebraically closed field, or where $\mathcal{G}$ is either split or anisotropic over $K$. 

Given a word $w=w(x_1,\ldots,x_k)\in F_d$ and a group $G$, the set of $w$-values is the union of the word map image with its inverse.

$$
w(G)=\{w(g_1,\ldots,g_k), w(g_1,\ldots,g_k)^{-1}\mid g_i\in G\}\tp
$$

The width of the word $w$ in the group $G$ is a minimal constant $C$ such that $\<w(G)\>=w(G)^C$, where $\<w(G)\>$ is the group generated by $w(G)$ and $w(G)^C$ is the set of elements that can be writen as a product of $C$ elements from $w(G)$.

To put the results of the present paper in a context, we now recall some results that are known at the moment.

$\bullet$ If $G$ is a connected simple algebraic group over an algebraically closed fied, and $w$ is a non-trivial word, then $w(G)$ contains an open dense subset, see~\cite{BorelFree}; hence the width of $w$ is at most 2.

$\bullet$ The set of word values in the unitary group $U_n$ can have arbitrarily small radii, see~\cite{Thom}. Hence the width of $w$ can be arbitrarily large (but it is finite for every $w$).

$\bullet$ For any word $w$, if a finite simple group $G$ is large enough (the meaning of "large enough" depends on $w$), then the width of $w$ in $G$ is at most 2, see~\cite{LarShal}.

$\bullet$ For any word $w$, if a classical connected real compact Lie group $G$ has large enough rank (the meaning of "large enough" depends on $w$), then the width of $w$ in $G$ is at most~2, see~\cite{HuiLarShal}.
 
$\bullet$ If $G$ is a Chevalley group over an infinite field, then the width of any word in $G$ is at most 5, also~\cite{HuiLarShal}.

$\bullet$ If $D$ is a skew-field that has finite dimension over its center, then the width of $w=[x,y]$ in $\GL(n,D)$ is at most [the width of $w=[x,y]$ in $D^*$ divided by $n$ (upper integer part)]+1 and at least $($[the width of $w=[x,y]$ in $D^*$]$+2n^2-3n+1):(8n^2-13n+8)$, see~\cite{GvozCommLength}, also see~\cite{KursovGL},\cite{KursovDiss},\cite{CommInGLD}, and \cite{EgorchDiss}.

$\bullet$ If $n$ is large enough, then the width of $w=[x,y]$ in $\SL(n,\Z)$ is at most 6, see~\cite{DenVasLength}. The question of whether it is equal to one is wide open.

$\bullet$ If $G(R)$ is a Chevalley group of rank at least 2 over a commutative ring $R$ of finite Bass--Serre dimension such that it coincides with its elementary subgroup, then the word $w=[x,y]$ has finite width if and only if $G(R)$ is boundedly generated by elementary root elements, see~\cite{HSVZwidth}.

$\bullet$ The width of any word in $\SL(n,\Z)$, $n\ge 3$ is bounded by some constant that depend only on $n$, see~\cite{AvniMeiri}.

$\bullet$ For any word $w$, if $n$ is large enough (the meaning of "large enough" depends on~$w$), then the width of $w$ in $\SL(n,\Z)$ is at most 87, also~\cite{AvniMeiri}.

\smallskip

In the present paper we study the width of words in Chevalley groups over a localisation of the ring of integers in a number field. Our main result is the following theorem, which is a generalisation of~\cite[Theorem~1.2]{AvniMeiri}  

{\thm\label{Main} Let $\Phi$ be an irreducible root system with $\rk\Phi\ge 2$. Let $d$ be a positive integer. Then there exists a constant $C=C(\Phi,d)$ such that for any ring $R$ satisfying:
\begin{enumerate}
\item $R$ is a localisation of the ring of integers in a number field $K$, with $[K\colon \Q]=d$;

\item if $\Phi=C_2$ of $G_2$, then $R$ has no residue field of two elements;

\end{enumerate}

the width of any word in $G_{\mathrm{sc}}(\Phi,R)$ does not exceed $C(\Phi,d)$.}

The paper is organised as follows. In Section~\ref{BasicNotationSec} we introduce the basic notation and give some necessary preliminaries. In Section~\ref{ProofSec} we prove Theorem~\ref{Main}.

I am grateful to Eugene Plotkin and Boris Kunyavskii for very useful discussions regarding various aspects of this work and permanent attention to it. 

\section{Basic notation}
\label{BasicNotationSec}

\subsection{Words and verbal width}

A word is an element in a free group. Given a word $w=w(x_1,\ldots,x_k)\in F_d$ and a group $G$, we have the word map $\map{w}{G^k}{G}$ defined by substitution. The set of $w$-values is the set

$$
w(G)=\{w(g_1,\ldots,g_k), w(g_1,\ldots,g_k)^{-1}\mid g_i\in G\}\tp
$$

The width of a word $w$ in a group $G$ is a minimal constant $C$ such that $\<w(G)\>=w(G)^C$. The verbal width of a group $G$ is the supremum of the widths of all words.

\subsection{Chevalley groups}

Let $\Phi$ be an irreducible root system in the sense of \cite{Bourbaki4-6}, $R$ a commutative associative ring with unity, $G(\Phi,R)=G_{\mathrm{sc}}(\Phi,R)$ a simply connected Chevalley group of type $\Phi$ over $R$, $T(\Phi,R)$ a split maximal torus of $G(\Phi,R)$. For every root $\alpha\in\Phi$ we denote by $X_\alpha=\{x_\alpha(\xi),\colon \xi\in R\}$ the corresponding unipotent root subgroup with respect to $T$. We denote by $E(\Phi,R)$ the elementary subgroup generated by all $X_\alpha$, $\alpha\in\Phi$. 

For an ideal $I\unlhd R$ we denote by $\rho_I$ the reduction homomorphism
$$
\map{\rho_I}{G(\Phi,R)}{G(\Phi,R/I)}\tc
$$
induced by projection $R\to R/I$.

The kernel of this homomorphism $G(\Phi,R,I)=\Ker(\rho_I)\le G(\Phi,R)$ is called the {\it principal congruence subgroup}. Set $E(\Phi,I)=\<x_\alpha(\xi)\tc\,\alpha\in\Phi\tc\,\xi\in I\>$. Then, by definition, the {\it relative elementary subgroup} $E(\Phi,R,I)$ of level $I$ is the normal closure of $E(\Phi,I)$ in $E(\Phi,R)$, actually it is automatically normal in $G(\Phi,R)$.

\subsection{Affine schemes and generic elements}

The functor $G(\Phi,-)$ from the category of rings to the category of groups is an affine group scheme (a Chevalley--Demazure scheme). This means that its composition with the forgetful functor to the category of sets is representable, i.e.,
$$
G(\Phi,R)=\Hom (\Z[G],R).
$$
The ring $\Z[G]$ here is called the {\it ring of regular functions} on the scheme $G(\Phi,-)$. 

An element $g_{\mathrm{gen}}\in G(\Phi,\Z[G])$ that corresponds to the identity ring ho\-mo\-mo\-rphism is called the {\it generic element} of the scheme $G(\Phi,-)$. This element has a universal property: for any ring $R$ and for any $g\in G(\Phi,R)$, there exists a unique ring homomorphism
$$
\map{f}{\Z[G]}{R},
$$
such that $f_*(g_{\mathrm{gen}})=g$. For details about application of the method of generic elements
to the problems similar to that of ours, see the paper of A.~V.~Ste\-pa\-nov~\cite{StepUniloc}.
 
\section{Proof of Theorem~\ref{Main}}
\label{ProofSec}

We need several lemmas. The next lemma is a simple technical group-theoretic result. 

{\lem\label{FiniteCover} Let $G$ be a group, and let $X\sub G$ be a symmetric set containing identity element such that there are $k$ translates of $X$ that cover $G$. Then we have $X^{3k-2}=\<X\>$.}
\begin{proof}
	Assume that $G=\bigcup_{i=1}^k g_iX$, without loss of generality we may assume that $g_1$ is the identity element. Let $g\in \<X\>$, and let $h_0$,$\ldots$,$ h_n$ be the shortest sequence such that $h_0$ is the identity element, $h_{j+1}\in h_jX$ for all $0\le j\le n-1$, and $h_n=g$. We must prove that $n\le 3k-2$.
	
	For all $0\le j\le n$ let $h_j\in g_{i_j}X$ with $i_0=1$. If for some $j$ and $j'$ we have $j'-j>2$ and $i_j=i_{j'}$, then our sequence can be shortened, because $h_{j'}\in h_jX^2$. In particular this means that for any $i$ we have $|\{j\mid i_j=i\}|\le 3$. Similarly, if for some $j>1$ we have $i_j=1$, then we can shorten our sequence by removing everything between $h_0$ and $h_j$; hence $|\{j\mid i_j=1\}|\le 2$.
	
	Therefore, we have $n+1\le 3(k-1)+2$, or $n\le 3k-2$.
\end{proof}

Let $G=G(\Phi,-)$ be a simply connected Chevalley--Demazure group scheme. Let $\Z[G]$ be the corresponding ring of regular functions. The center of $G(\Phi,-)$ is a closed subscheme, so it corresponds to some ideal in $\Z[G]$. Let us fix some generators $\{\eta_i\}_{i=1}^{i_{\max}}$ of that ideal. For example we can take $\eta_i$ to be the differences between entries of the adjoint representation with the corresponding entries of the identity matrix.

Further we need an effective version of the well known result by Abe about classification of normal subgroups in Chevalley groups (see \cite{AbeNormal}).

{\lem \label{EfficientSandwich} Let $\Phi$ be an irreducible root system with $\rk\Phi\ge 2$. Set $e=1$ if $\Phi$ is of type ADE, $e=2$ if $\Phi$ is of type BCF, and $e=3$ if $\Phi=G_2$. Let $\eta_i$ be as above. Then the following holds.
	
\begin{enumerate} \item Let $\Phi$ be distinct from $C_2$ and $G_2$. Then there exists a constant $C=C(\Phi)$ such that for any commutative ring $R$ and for any $g\in G(\Phi,R)$ each of the following elements can be decomposed as a product of at most $C$ elements conjugate to $g$ or~$g^{-1}$:
\begin{itemize}
\item elements $x_{\alpha}(\xi\eta_i(g))$, $\xi\in R$, $\alpha\in \Phi$ if $\Phi$ is of type ADE;

\item elements $x_{\alpha}(\xi\eta_i(g))$, $\xi\in R$, $\alpha$ is a short root if $\Phi$ is not of type ADE;

\item elements $x_{\alpha}(\xi\eta_i(g)^e)$ and $x_{\alpha}(e\xi\eta_i(g))$, $\xi\in R$, $\alpha$ is a long root if $\Phi$ is not of type ADE.
\end{itemize} 

\item Let $\Phi=C_2$,$G_2$, and let $l$ be a positive integer. Then there exists a constant $C=C(\Phi,l)$ such that for any commutative ring $R$ containing elements $t_1$,$\ldots$,$t_l$ such that the set $\{t_i^2+t_i\colon 1\le i\le l\}$ generates a unit ideal and for any $g\in G(\Phi,R)$ each of the elements as above can be decomposed as a product of at most $C$ elements conjugate to $g$ or $g^{-1}$.
\end{enumerate}}
\begin{proof}
	First let us prove that the elements in question must belong to the normal subgroup generated by $g$.
	
	Let $H\unlhd G(\Phi,R)$ be the normal subgroup generated by $g$. By \cite[Theorem~3]{AbeNormal} there exists a so called admissible pair $(A,B)$ such that
	$$
	E(R,A,B)\le H\le E^*(R,A,B)\tc
	$$ 
	where $E(R,A,B)$ is a normal subgroup of $E(\Phi,R)$ generated by $x_{\alpha}(\xi)$, $\xi\in A$, $\alpha$ is a short root, and $x_{\beta}(\zeta)$, $\zeta\in B$, $\beta$ is a long root (in case where $\Phi$ is of type ADE, we have $A=B$); and 
	$$
	E^*(R,A,B)=\{h\in G(\Phi,R)\mid [h, E(\Phi,R)]\le E(R,A,B)\}\tp
	$$

	In case $\Phi=C_2$,$G_2$ this theorem requires that $R$ has no residue field of two elements, which we ensure by requiring that the set $\{t_i^2+t_i\colon 1\le i\le l\}$ generate a unit ideal.

	By definition of an admissible pair $A$ is an ideal, and $B\sub A$; hence $E(R,A,B)$ is trivial modulo $A$. Therefore, modulo $A$, the element $g$ commutes with the elementary subgroup, and hence must be central (see, for example \cite[Main Theorem]{AbeHurley}). Therefore, we have $\eta_i(g)\in A$ for all $g$, which proves that the elements of the first and the second type are in~$H$.
	
	Also by definition of an admissible pair, the ideal generated by $e\xi$ and $\xi^e$ for $\xi\in A$ is contained in $B$, which proves that the elements of the third type are in $H$.
	
	Now let us prove that the number of conjugates of $g$ and $g^{-1}$ required to express these elements can be bounded by a constant that depends only on $\Phi$ or, if $\Phi= C_2$,$G_2$, on $\Phi$ and $l$.
	
	We start with the case $\Phi\ne C_2$,$G_2$. Consider the ring $\Z[G][\xi]$, i.e the ring of polynomials over the ring of regular functions. Let $g=g_{\mathrm{gen}}$ be the generic element. By what we proved previously, the elements $x_{\alpha}(\xi\eta_i)$ (or $x_{\alpha}(\xi\eta_i^e)$ and $x_{\alpha}(e\xi\eta_i)$ in case of long roots in multiply laced systems) are products of conjugates of $g$ and $g^{-1}$. Take $C$ to be the maximum of the corresponding numbers of conjugates. Since $g$ can be specialised to any element of the group over any ring and simultaneously the variable $\xi$ can be specialised to any element of this ring, it follows that the bound $C$ works in all the cases.
	
	In case  $\Phi=C_2$,$G_2$ the proof is similar. The only difference is that we must consider the ring $\Z[G][\xi,\zeta_1,\ldots,\zeta_k,t_1,\ldots,t_l]/(\zeta_1(t_1^2+t_1)+\ldots+\zeta_k(t_k^2+t_l)-1)$.
\end{proof}

The next lemma helps to use the previous one in the cases $\Phi=C_2$,$G_2$.

{\lem\label{NoF2} Let $R$ be a localisation of the ring of integers in a number field, and suppose that $R$ has no residue field of two elements. Then there are $t_1$,$t_2\in R$ such that $t_1^2+t_1$ and $t_2^2+t_2$ generate the unit ideal.}
\begin{proof}
	Take $t_1$ to be any element such that $t_1^2+t_1\ne 0$. Let $J$ be the Jacobson radical of the ideal generated by $t_1^2+t_1$. Then $R/J$ is a product of fields, none of which is $\F_2$. Hence we can take $t_2$ to be the element such that the images of $t_2^2+t_2$ in all these fields are non zero. Therefore, $t_1^2+t_1$ and $t_2^2+t_2$ generate the unit ideal.
\end{proof}
 
 For a ring $R$ and an ideal $I\unlhd R$ the group $E(\Phi,R,I)$ is generated by elements $x_{\alpha}(\xi)^g$, where $\alpha\in\Phi$, $\xi\in I$ and $g\in G(\Phi,R)$. We denote by $E(\Phi,R,I)^{\le C}$ the subset of $E(\Phi,R,I)$ that consist of elements that can be decomposed as a product of at most $C$ such generators.
 
 The next two lemmas allow one to establish the equality $E(\Phi,R,I)=E(\Phi,R,I)^{\le C}$ for arithmetic rings. See also \cite{GvozWidth} and \cite{VavBoundedRelative} for some different variations of it.
 
 {\lem\label{Logic} Let $\Phi$ be an irreducible root system with $\rk\Phi\ge 2$ and $d$ be a positive integer. Then there exists a first-order statement $\ph=\ph_{\Phi,d}$ in the language of rings with the following properties:
 \begin{enumerate}
 \item If $R$ is a localisation of the ring of integers in a number field $K$, with $[K\colon \Q]=d$, then $\ph$ holds in $R$.
 
 \item If the statement $\ph$ holds in a ring $R$ and $I\unlhd R$ is a 2-generated ideal, then $[G(\Phi,R,I):E(\Phi,R,I)]\le 4(8d)!$.
\end{enumerate}}
 
 \begin{proof}
 	Pick a long root $\alpha\in\Phi$. By~\cite[Corollary 4.5]{Matsumoto69}, if $R$ is a localisation of the ring of integers in a number field, then we have $G(\Phi,R,I)=G(\{\alpha,-\alpha\},R,I)E(\Phi,R,I)$. The proof in~\cite{Matsumoto69} is constructive, so it is possible to extract from this proof a bound $C=C(\Phi)$ such that $G(\Phi,R,I)=G(\{\alpha,-\alpha\},R,I)E(\Phi,R,I)^{\le C}$.
 	
 	Note that for fixed $\Phi$ and $C$ the statement "for any 2-generated ideal $I$, we have $G(\Phi,R,I)=G(\{\alpha,-\alpha\},R,I)E(\Phi,R,I)^{\le C}$" is a first-order statement in the language of rings. Denote this statement by $\mathrm{Red}(\Phi,C)$.
 	
 	Now we construct our statement $\ph$ as a conjunction of $\mathrm{Red}(\Phi,C)$ and the statements $\mathrm{SR}_{1\tfrac{1}{2}}$, $\mathrm{Gen}(2,1)$, $\mathrm{Gen}(2(8d)!,1)$, $\mathrm{Exp}(2(8d)!,2)$ (see~\cite[Definition~2.10, 3.2, 3.6]{Morris}).
 	
 	If $R$ is a localisation of the ring of integers in a number field $K$, with $[K\colon \Q]=d$, then it satisfies $\mathrm{Red}(\Phi,C)$ as explained above, and it satisfies the other statements by \cite[Lemma~2.13, Corrollary~3.5, Theorem~3.9]{Morris}.
 	
 	Now let $R$ be a ring that satisfies $\ph$. Since $\rk\Phi>1$, we have one of the following.
 	
 	{\bf Case 1} The root $\alpha$ is contained in a subsystem of type $A_2$. Then it follows from $\mathrm{Red}(\Phi,C)$ that the quotient $G(\Phi,R,I)/E(\Phi,R,I)$ is generated by relative Mennicke symbols; and by \cite[Theorem 3.11]{Morris} the order of the universal Mennicke group is bounded by $2(8d)!$. 
 	
 	{\bf Case 2} The root $\alpha$ is contained in a subsystem of type $C_2$. Then it follows from $\mathrm{Red}(\Phi,C)$ that we have a surjection  $G(C_2,R,I)/E(C_2,R,I)\to  G(\Phi,R,I)/E(\Phi,R,I)$ (in other words $\mathrm{Red}(\Phi,C)$ implies the surjective stability of the $K_1$-functor); and by \cite[Theorem 3.16]{TrostChevalley}  the order of the group $G(C_2,R,I)/E(C_2,R,I)$ is bounded by $4(8d)!$. 
 \end{proof}

{\lem\label{WidthOfERI} Let $\Phi$ be an irreducible root system with $\rk\Phi\ge 2$ and $d$ be a positive integer. Then there exists a constant $C=C(\Phi,d)$ such that for any ring $R$ that is a localisation of the ring of integers in a number field $K$, with $[K\colon \Q]=d$, and any ideal $I\unlhd R$ we have $E(\Phi,R,I)=E(\Phi,R,I)^{\le C}$.}

\begin{proof}
	
By Lemma~\ref{FiniteCover} it is enough to show that there exists a constant $C(\Phi,d)$ such that for any ring $R$ as above and any ideal $I\unlhd R$ there exist $g_1$,$\ldots$,$g_{4(8d)!}$ such that $G(\Phi,R,I)=\bigcup_i g_iE(\Phi,R,I)^{\le C}$.

Assume that the above statement is false. Then for every $C\in \N$ there exist a ring~$R_C$ as above, an ideal $I_C\unlhd R_C$ and elements $g_{C,1}$,$\ldots$,$g_{C,4(8d)!+1}\in G(\Phi,R_C,I_C)$ such that $(g_{C,i})^{-1}g_{C,j}\notin E(\Phi,R_C,I_C)^{\le C}$ if $i\ne j$. Here the elements $g_{C,i}$ are constructed as follows: $g_{C,1}$ is picked arbitrary and $g_{C,j+1}$ is picked from $G(\Phi,R_C,I_C)\sm \bigcup_{i\le j} g_iE(\Phi,R,I)^{\le C}$ that guarantees that  $(g_{C,i})^{-1}g_{C,j}\notin E(\Phi,R_C,I_C)^{\le C}$ for $i<j$, and since $E(\Phi,R_C,I_C)^{\le C}$ coincides with its inverse, we have the same for  $i>j$.   

It is well known fact from model theory that the first order statements are inherited by ultraproducts. Choose a non-principal ultrafilter $\mathcal{U}$ on $\N$, and let $R_{\infty}$ be the ultraproduct of $R_C$ over~$\mathcal{U}$. Then the statement $\ph$ from Lemma~\ref{Logic} holds in $R_{\infty}$ and $G(\Phi,R_{\infty})$ is isomorphic to the ultraproduct of $G(\Phi,R_C)$ over $\mathcal{U}$. Let $I_{\infty}$ be the ideal of $R_{\infty}$ represented by $\prod_C I_C$; this ideal is 2-generated because all the $I_C$ are 2-generated in $R_C$. For every $1\le i\le 4(8d!)+1$, let $g_i\in G(\Phi,R_{\infty},I_{\infty})$ be the element represented by $(g_{C,i})_C$. Then $g_1$,$\ldots$,$g_{(4(8d!)+1}$ belong to different cosets of $E(\Phi,R_{\infty},I_{\infty})$ contradicting Lemma~\ref{Logic}. 
\end{proof}

Now we combine the previous results to establish the following lemma.

{\lem\label{ExistsAnIdeal} Let $\Phi$ be an irreducible root system with $\rk\Phi\ge 2$. Let $d$ be a positive integer. Then there exists a constant $C=C(\Phi,d)$ such that for any ring $R$ satisfying assumptions (1)--(2) of Theorem~\ref{Main} and any symmetric conjugation-invariant subset $X\sub G_{\mathrm{sc}}(\Phi,R)$ that contains the identity element, and that is not contained in the center there is a non-trivial ideal $I\unlhd R$ such that $E(\Phi,R,I)\le X^C$.}
\begin{proof}
Let $C_1$ be the constant from Lemma~\ref{EfficientSandwich} (in case $\Phi=C_2$,$G_2$ we take $l=2$), and let $C_2$ be the constant from Lemma~\ref{WidthOfERI}. Then we can take $C=C_1C_2$.

Indeed, let $X$ be as above. Since $X$ is not contained in the center, it follows that there are $x\in X$ and $i$ such that $\eta_i(x)\ne 0$. Take $I=e\eta_i(x)R$, then by Lemmas~\ref{EfficientSandwich} and~\ref{NoF2} all the elements $x_{\alpha}(\xi)$, $\alpha\in\Phi$, $\xi\in I$ are in $X^{C_1}$. Therefore, $E(\Phi,R,I)\le X^{C_1C_2}$.

\end{proof}

For a ring $R$ and a positive integer $n$ we denote by $\mu_n(R)$ the group of $n$-th roots of unity in the ring $R$.

For a ring $R$ and an ideal $\P\unlhd R$ we denote by $R_{(\P)}$ the corresponding completion of $R$. Note that $R_{(\P)}$ is a DVR.

Now we prove some technical lemmas about the groups over such completions.

{\lem\label{RootsOfUnity}  Let $R$ be a localisation of the ring of integers in a number field $K$, with $[K\colon \Q]=d$, and let $n$ be a positive integer. Then for any prime ideal $\P\unlhd R$ and any positive integer $m$ we have $|\mu_n(R/\P^m)|\le n^{2d}$.}
\begin{proof}
	
	Let $p$ be the characteristic of $R/\P$, let $q=|R/\P|$. Let $n=p^kz$, where $z$ is coprime to $p$. Note that if $k=0$, then by Hensel's lemma every $n$-th root of unity in $\mu_n(R/\P^m)$ can be lifted to $R_{(\P)}$ and hence $|\mu_n(R/\P^m)|\le n$; thus we may assume that $k\ge 1$. Let $\pi$ be the uniformizing element in $R_{(\P)}$, and let $p=\xi \pi^r$, where $\xi\in R_{(\P)}^*$.
	
	We claim that if $m>2kr$, then we have
	$$
	\Ker\left((R/\P^{m+1})^*\to (R/\P^{m})^* \right)\le ((R/\P^{m+1})^*)^n\tp
	$$
	Indeed, elements from this kernel are of the from $1+s\pi^m$. We have
	$$
	\left(1+\frac{s}{z\xi^k}\pi^{m-kr}\right)^n\equiv 1+s\pi^m \mod \P^{2(m-kr)}\tc
	$$
	and if $m>2kr$, then $2(m-kr)>m$.
	 
	Therefore, for $m>2kr$, we have 
	$$
	(R/\P^{m+1})^*/((R/\P^{m+1})^*)^n\simeq (R/\P^{m})^*/((R/\P^{m})^*)^n\tp
	$$
	
	Since $(R/\P^{m})^*$ is a finite abelian group, we have
	$$
	|\mu_n(R/\P^{m})|=|(R/\P^{m})^*/((R/\P^{m})^*)^n|\tp
	$$
	
	Therefore, for $m>2kr$, we have 
	$$
	|\mu_n(R/\P^{m+1})|=|\mu_n(R/\P^{m})|\tp
	$$
	
	Since $(R/\P)^*$ is a cyclic group of order coprime to $p$, we have $|\mu_n(R/\P)|\le z$. For the kernel of the reduction map we have: 
	
	$$
	|\Ker(\mu_n(R/\P^{m})\to \mu_n(R/\P))|\le |\Ker(R/\P^{m}\to R/\P)|=q^{m-1}\tp
	$$ 
	
	Hence we have $|\mu_n(R/\P^{m})|\le zq^{m-1}$. Therefore, for every $m$ we have
	
	$$
	|\mu_n(R/\P^{m})|\le zq^{2kr}\le zp^{2kd}\le (zp^k)^{2d}=n^{2d}\tp
	$$
	
	Here we need to explain the equality $q^r\le p^d$. Let $L$ be the fraction field of $R_{(\P)}$. Let $F$ be the maximal unramified extension of $\Q_p$ in $L$. Then $[L:\Q_p]=[L:F][F:\Q_p]=r\cdot [\F_q:\F_p]$. Therefore, $q^r=p^{[L:\Q_p]}$. On the other hand $d$ is a rank of $R$ as a $\Z$-module, hence it is a rank of $R\otimes_{\Z}\Z_{p}$ as a $\Z_p$-module. Now $R\otimes_{\Z}\Z_{p}=\prod_{\P_i} R_{(\P_i)}$ where $\P_i$ are all primes of $R$ lying over $p$. Therefore, we have $d=\sum_{\P_i} \rk_{\Z_p} R_{(\P_i)}=\sum_{\P_i} [\mathrm{Frac}(R_{(\P_i)}):\Q_p]$, and $[L:\Q_p]$ here is one of the summands.
\end{proof}

{\lem\label{InCompletion} Let $\Phi$ be an irreducible root system with $\rk\Phi\ge 2$. Let $d$ be a positive integer. Then there exists a constant $C=C(\Phi,d)$ such that for any ring $R$ satisfying assumptions (1)--(2) of Theorem~\ref{Main}, any prime ideal $\P\unlhd R$ and any symmetric conjugation-invariant subset $X\sub G_{\mathrm{sc}}(\Phi,R_{(\P)})$ that contains the identity element we have $\<X\>=X^C$.}

\begin{proof}
Let $C_1$ be the constant from Lemma~\ref{EfficientSandwich} (in case $\Phi=C_2$,$G_2$ take $l=1$; note that in a local ring with the residue field distinct from $\F_2$ there exists an element $t$ such that $t^2+t$ is invertible). Let $n$ be a number such that the center of the group scheme $G(\Phi,-)$ is the group $\mu_n$ of the $n$-th roots of unity (in case $\Phi=D_{2l}$, where the center is $\mu_2\times\mu_2$ we take $n=4$), so $n$ depends only on $\Phi$.

We claim that we can take
$$
C=C_1\cdot\left(\frac{3}{2}|\Phi|+2\rk\Phi+1\right)\cdot\left( 3^{d(|\Phi|+\rk\Phi)+1}\cdot n^{2d}-2\right)\tp
$$

Indeed, let $X$ be as above. We may assume that $X$ is non-central, because otherwise, we have $|\<X\>|\le n$ and then, clearly, $\<X\>=X^n$. Let $k$ be the biggest number such that $X\sub CG(\Phi,R_{(\P)},\P^k)$, where $CG(\Phi,R_{(\P)},\P^k)$ is the full congruence subgroup, i.e. the preimage of the center under the reduction map $G(\Phi,R_{(\P)})\to G(\Phi,R_{(\P)}/\P^k)$.

Recall that in the notation of Lemma~\ref{EfficientSandwich} we set $e=1$ if $\Phi$ is of type ADE, $e=2$ if $\Phi$ is of type BCF, and $e=3$ if $\Phi=G_2$.

{\bf Case 1} $k=0$.

Recall that by $\eta_i$ we denote some fixed generator of the ideal in $\Z[G]$ that correspond the scheme center of $G$. In this case there exists $x\in X$ and an index $i$ such that $\eta_i(x)\notin\P$. Thus $\eta_i(g)^e$ is invertible; hence $x_{\alpha}(\xi)\in X^{C_1}$ for all $\alpha\in\Phi$, $\xi\in R_{(\P)}$, and by \cite[Lemma~5.3]{SinchukSmolensky} (take $I$ to be the unit ideal) we have

$$
G(\Phi,R_{(\P)})=X^{C_1\left(\frac{3}{2}|\Phi|+2\rk\Phi+1\right)}\tp
$$

{\bf Case 2} $k\ge 1$. 

Let $\pi$ be the uniformizing element, $p$ be the characteristic of $R/\P$, let $q=|R/\P|$, and let $p=\xi \pi^r$, where $\xi\in R_{(\P)}^*$.

There is an element $x\in X$ and an index $i$ such that $\eta_i(x)\notin\P^{k+1}$. Hence $x_{\alpha}(\xi e \pi^k)\in X^{C_1}$ for all $\alpha\in\Phi$, $\xi\in R_{(\P)}$, and by \cite[Lemma~5.3]{SinchukSmolensky} (take $I=e\P^{k}$) we have
$$
G(\Phi,R_{(\P)},e\P^{k})\le X^{C_1\left(\frac{3}{2}|\Phi|+2\rk\Phi+1\right)}\tp
$$

Therefore, since $X\sub CG(\Phi,R_{(\P)},\P^k)$, it follows by Lemma~\ref{FiniteCover} that it is enough to show that $[CG(\Phi,R_{(\P)},\P^k):G(\Phi,R_{(\P)},e\P^{k})]\le 3^{d(|\Phi|+\rk\Phi)}\cdot n^{2d}$.

It follows by Lemma~\ref{RootsOfUnity} that $[CG(\Phi,R_{(\P)},\P^k):G(\Phi,R_{(\P)},\P^{k})]\le n^{2d}$. Therefore, it remains to estimate $[G(\Phi,R_{(\P)},\P^k):G(\Phi,R_{(\P)},e\P^{k})]$. From now we may assume that $p=e\in\{2,3\}$ because otherwise this index is equal to one.

All elements from $G(\Phi,R_{(\P)},\P^k)$ belong to the main Gauss cell, which is an open subscheme of $G(\Phi,-)$ isomorphic to $\mathbb{A}_1^{|\Phi|}\times\mathbb{G}_m^{\rk\Phi}$. Modulo $e\P^{k}=\P^{k+r}$ every component of this product can take one of the $q^r$ values; hence
$$
[G(\Phi,R_{(\P)},\P^k):G(\Phi,R_{(\P)},e\P^{k})]\le q^{r(|\Phi|+\rk\Phi)}\le p^{d(|\Phi|+\rk\Phi)}\le 3^{d(|\Phi|+\rk\Phi)}\tp
$$
\end{proof}

We need one more technical group-theoretic lemma.

{\lem\label{Product} Let $G_i$ be a collection of groups, and let $G=\prod_i G_i$. Let $w$ be a word. Assume that for all $i$ the width of $w$ in $G_i$ is at most $C$. Then the width of $w$ in $G$ is at most $2C$.}

\begin{proof}
	Let $g\in \<w(G)\>$, and let $g_i$ be the components of $g$. Set
	$$
	w^+(G)=\{w(h_1,\ldots,h_k)\mid h_i\in G\}\tc\qquad 	w^-(G)=\{w(h_1,\ldots,h_k)^{-1}\mid h_i\in G\}\tp
	$$
	
	By assumption each of the $g_i$ is a product of at most $C$ elements from $w(G_i)$. Since $w^+(G_i)$ is conjugation-invariant, we can reorder this product so that the elements from $w^+(G_i)$ comes before the elements of $w^-(G_i)$. Therefore, $g_i\in w^+(G_i)^C w^-(G_i)^C$; hence $g\in w^+(G)^C w^-(G)^C\sub w(G)^{2C}$.
\end{proof}

Now let us prove Theorem~\ref{Main}.

\smallskip

Let $X=w(G(\Phi,R))$. Clearly, the set $X$ is not contained in the center (because $G(\Phi,R)$ contains a copy of $\SL(2,\Z)$, and hence a copy of a free subgroup). Let $R$ be the ring that satisfy (1)--(2). Let $C_1$ be the constant from Lemma~\ref{ExistsAnIdeal}, $C_2$ be the constant from Lemma~\ref{InCompletion}, and $m$ be the number of roots of unity in $K$. We prove that $\<X\>=X^{(C_1+2C_2)(3m-2)}$. Clearly $m$ is bounded by some constant that depens on $d$ (we have $\ph(m)\le d$, where $\ph$ is the Euler function; hence $m\le 2d^2$); therefore, that would be enough.

By Lemma~\ref{ExistsAnIdeal}, there is a non-trivial ideal $I\unlhd R$ such that $E(\Phi,R,I)\le X^{C_1}$. By \cite[Corollary 4.3]{BassMilnorSerre} and \cite[Corollary 4.5]{Matsumoto69} we have $|G(\Phi,R,I)/E(\Phi,R,I)|\le m$. Therefore, by Lemma~\ref{FiniteCover}, it is enough to prove that $\<X\>\le G(\Phi,R,I)X^{2C_2}$.

Let $I=\prod_{i=1}^{n} \P_i^{k_i}$, and let $g\in \<X\>$. It follows by Lemmas~\ref{InCompletion} and~\ref{Product} the image of~$g$ in $G(\Phi,\prod_{i=1}^{n} R_{(\P_i)})=\prod_{i=1}^{n} G(\Phi,R_{(\P_i)})$ belong to $w(G(\Phi,\prod_{i=1}^{n} R_{(\P_i)}))^{2C_2}$. Since $R/I\simeq (\prod_{i=1}^{n} R_{(\P_i)})/\prod_{i=1}^{n} \P_i^{k_i}$, and the reduction map $G(\Phi,R)\to G(\Phi, R/I)$ is surjective, it follows that $g$ is congruent to an element from $X^{2C_2}$ modulo $G(\Phi,R,I)$ as required.


\end{document}